\documentclass[11pt]{article}

\usepackage{euscript,amsmath,amssymb ,graphicx,setspace}
\usepackage{srcltx}

\title{On the morphisms of fractal curves that increase their smoothness.}

\author{Kirill Kamalutdinov,  Svetlana Sorokina}

\newcommand {\Ga} {\Gamma}

\newcommand{\e}{\varepsilon} 

\newcommand{\eS}{{\EuScript S}}

\newcommand{\eT}{{\EuScript T}}

\newcommand{\eZ}{{\EuScript Z}}

\newcommand{\rr}{\mathbb{R}}

\begin{document}

\maketitle

\begin{abstract} 
We propose a construction  which  transforms   a  self-similar zipper in $\mathbb R^n$  to  a  self-affine  zipper $\mathbb R^{n+1}$  whose  attractor  is  a  smooth curve.
\end{abstract}

\vspace{12mm}\\

Smooth  self-affine curves may  be  constructed in two ways:

1. As it was proved  by M.Barnsley [3],  for some  special types  of fractal interpolation  functions  their  integration  gives  differentiable   functions  whose  graphs  are  self-affine;\\
2. The other  way  is  to use the  algorythm of building self-affine curves proposed  by Alexey Kravchenko in 2005. [5]

We propose  one  more construction,  which  transforms   a  self-similar zipper  to  a  self-affine  zipper whose  attractor  is  a  smooth curve.

\bigskip

{\bf Definition 1.} {\it Let X be  a  complete  metric  space. A system $\eS=\{S_1,...,S_m\}$ of contraction mappings of $X$  to itself is called a {\em zipper} with vertices $\{z_{0},...,z_{m}\} $ and signature $\e=(\e_{1},..,\e_{m}), \e_i\in\{0,1\}$, if for  any $i=1,...,m$,   $S_i(z_{0})=z_{i-1+\e_{i}}$ and $S_{i}(z_{m})=z_{i-\e_{i}}$. } [3]

The zipper $\eS$   is called  self-similar or self-affine, if all $S_i$ are similarities or affine mappings.

{\bf Definition 2.} {\it A compact set $K \subset X$ is called an attractor or invariant set of the system $\eS$  if $K=\bigcup\limits_{i=1}^{n} s_i (K)$. }[3]

An attractor of any system exists and unique due to Hutchinson theorem [4]. Attractor of any zipper is arcwise connected and locally arcwise connected [1].

{\bf Definition 3.} {\it A zipper $\eS$ in $[0,1]$ with vertices $\{0=t_{0},...,t_{m}=1\}$,   is called a line zipper. }

 Suppose $\eS$ is a zipper with vertices $\{z_{0},...,z_{m}\}$ and   signature $\e$ and $\gamma$ is its attractor. Let $\eT=\{T_1,...,T_m\}$ be a line zipper with vertices $\{t_{0},...,t_{m}\}$ and  signature $\e$. 
 
 As it  was  proved  by Aseev et al [1], there is unique continuous $f: [0,1]\to \gamma$ such that $f(t_i)=z_i$ ($i=1,...,m$), and  for  any $i=1,...,m$  and $t\in[0,1]$, $f(T_i(t))=S_i(f(t))$.

    Such map $f$  is called a {\it linear parametrization} of the zipper $\eS$ by the zipper $\eT$.

{\bf Proposition 1.} {\it 
Let $\eS$ be a self-similar zipper in $\rr^{n}$ with vertices $\{z_{0},...,z_{m}\}$ and signature $\e$, $f$ be its linear parametrization by a line zipper $\eT$ with vertices $\{t_{0},...,t_{m}\}$.
Then $\eZ=\{S_1\times T_1,...,S_m\times T_m\}$ is a self-affine zipper in $\mathbb{R}^{n+1}$ with vertices ${(z_{0},t_0),...,(z_{m},t_m)}$ and signature $\e$, and its attractor is the graph $\Ga=\{(t,f(t)), t\in[0,1]\}$ of the function $f$.}   $\blacksquare$

From the construction of  the  function $f$ it  follows  that for  any $i=1,..,m$ and   any $t\in [t_{i-1},t_i]$
$$f(t)=S_i(f(T_i^{-1}(t))  \mbox{    where   }  T_i:[0,1]\to [t_{i-1},t_i]  \eqno{(1)} $$

Thus, $f$ can be  considered  as  a  fractal interpolation  function having  variable  signature $\e$. [3,4]

Taking  integral of  $f$ with respect  to $t$,  we  obtain a differentiable  function $g(t)=\int\limits_0^t f(\tau)d\tau$ on $[0,1]$,
whose  graph $\hat\Ga=\{t,g(t)\}$ is  a  self-affine Jordan  arc.

{\bf Proposition 2.} {\it Let $\eS$ be a self-similar zipper in $\mathbb{R}^{n}$ with vertices ${z_{0}=0,...,z_{m}}$ and signature $\e$, and $f$ be its linear parametrization by line zipper $\eT$ with vertices ${t_{0},...,t_{n}}$.
Then the graph of the function $g(t)=\int \limits_{0}^{t} f(\tau)d\tau$ is the attractor of a self-affine zipper in $\mathbb{R}^{n+1}$ with the signature $\e$.}

{\bf Proof.}
Put $z_m=b,\ g(1)=h$.

We write $S_i (z)= \begin{cases} z_{i-1} +A_i z, \mbox{  if  } \e_i=0 \\ z_{i} -A_i z, \mbox{     if  }  \e_i=1 \end{cases}$.

Here $A_i$ are the  similarities  sending $ 0,  b $ to $ 0$ and $z_i-z_{i-1}$ respectively. Denote  by $q_i$ the  contraction ratio $t_i-t_{i-1}$ of $T_i$. Then the  function $f$ satisfies  the  equations:
$$f(t)= \begin{cases} z_{i-1} + A_i f \left((t-t_{i-1})/q_i\right), \mbox{     if  }  \e_i=0 \\ 
z_{i-1} + A_i (b-f \left((t_{i}-t)/q_i\right),  \mbox{     if  }  \e_i=1 \end{cases} \mbox{    for   }t\in[t_{i-1},t_i]    \eqno{(2)}$$

  Integrating  from $t_{i-1}$ to $t$, we obtain:

$$g(t)-g(t_{i-1})=\begin{cases} z_{i-1}(t-t_{i-1}) +q_iA_i g\left( (t-t_{i-1})/q_i \right) \mbox{   if   }\e_i=0 ,\\  
z_{i}(t-t_{i-1}) +q_i A_i  \left( g\left((t_{i}-t)/q_i \right) - h \right) \mbox{   if   }e_i=1 \end{cases}  \eqno{(3)}$$

 $$g(t_i)-g(t_{i-1})=\begin{cases} z_{i-1}q_i +q_iA_i h, \e_i=0 \\ z_{i}q_i -q_iA_i h, \e_i=1 \end{cases}   \eqno{(4)}$$

Taking the  sum for all $i$, we see that $h$ satisfies:
$$\left(Id-\sum\limits_{i=1}^{m} (-1)^{\e_i}q_i A_i\right) h=\sum\limits_{i=1}^{m} q_i z_{i-1+\e_i}  \eqno{(5)}$$

 Finding $h$ from this  equation we find  the  values  of  all $g(t_i)$. 

From the  equation 3 we  see 
 that the graph $\tilde\Ga$ of   the  map $g$ is the attractor of $\eZ=\{W_1,...,W_m\}$, where
 
 $W_i (\tilde z)= \left( \begin{array}{cc} t_{i-1+\e_i} \\ g(t_{i-1+\e_i}) \\ \end{array} \right) + (-1)^{\e_i}q_i \left( \begin{array}{cccc} 1 & 0 & ... & 0 \\ z_{i-1+\e_i,\ 1} & & & \\ ... & & (-1)^{\e_i}A_i & \\ z_{i-1+\e_i,\ n} & & & \\ \end{array} \right) \tilde z $
  
  This   shows  that $\eZ$ is  a  self-affine  zipper  with vertices $(t_i,g(t_i))$ and  signature $\e$.
$\blacksquare$

{\bf Proposition 3.} { \it The  projection of $\tilde\Ga$ to $\rr^n$ is  smooth  at any  point, except 0.}

{\bf Proof.}
The projection of $\tilde\Ga$ to $\rr^n$ is $\{g(t):\ t\in[0,1]\}$. Taking $g'(t)=f(t)$ we obtain $g'(t)=0$ only when $t=0$. $\blacksquare$
\\

In this paper we consider the simplest zippers consisting of two mappings.

{\bf Example 1.} First, consider the line zipper, consisting of the similarities $S_{1}: [0,1] \rightarrow [0,p]$  and $S_{2}: [0,1] \rightarrow [p,1]$ , where $0 < p < 1$. Obviously, its attractor is the segment $[0,1]$. Parametrize it so that $f\left(\dfrac{1}{2}\right) = p$. The graph of this parameterization is the attractor of the zipper on the plane, consisting of two mappings that transform a single square in the rectangles $\left[0,\dfrac{1}{2}\right]\times[0,p]$ and $\left[\dfrac{1}{2},1\right]\times[p,1]$ (Fig. 1).

\begin{center}
\includegraphics {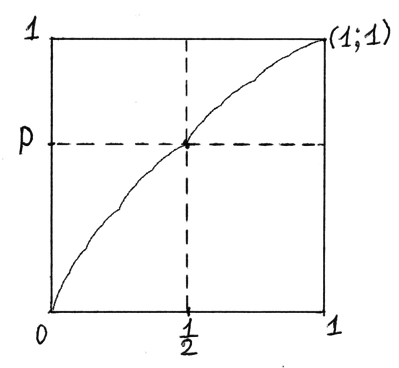}

Figure 1.
\end{center}

The parameterization of zipper satisfies equation:

\begin {center}
$f(x) =\begin {cases} pf(2x) , x\in \left[0;\dfrac{1}{2}\right]\\
(1-p)f(2x-1) + p , x\in \left[\dfrac{1}{2};1\right] \end {cases}$
\end {center}

Let $g(x) = \int\limits_{0}^{x} f(\xi)d\xi$ .

Then $g(x)$ satisfies equation:

\begin {center}
$g(x) =
\begin {cases}
\dfrac{p}{2}g(2x) , x\in \left[0;\dfrac{1}{2}\right]\\
\dfrac{(1-p)}{2}g(2x-1) + pt + \dfrac{p(p-1)}{2} , x\in \left[\dfrac{1}{2};1\right]
\end {cases}$

\end {center}

The zipper that specifies the graph (Fig. 2) of the function $g(x)$ consists of :

\begin {center}
$\begin {cases}
W_{1}(x,y) = \left( \dfrac{x}{2}; \dfrac{p}{2}y\right)\\
W_{2}(x,y) = \left( \dfrac{x}{2} + \dfrac{1}{2}; \dfrac{1-p}{2}y + \dfrac{p}{2}x + \dfrac{p^2}{2} \right)
\end {cases}$
\end {center}

\begin{center}
\includegraphics[scale=5]{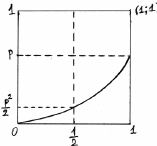}

Figure 2.
\end{center}

Put $y_i=f(x_i)$, $p_i=\dfrac{y_i-y_{i-1}}{y_m}$,  $y_0=0$, $q_i=x_i-x_{i-1}$ for $i=0,1,2$. Then

$f(x)=\begin{cases} p_1 f(\dfrac{x}{q_1}) \\ p_2 f(\dfrac{x-x_{1}}{q_2}) + y_{1} \end{cases}$, 
$g(x)=\begin{cases} p_1 q_1 g(\dfrac{x}{q_1}) \\ p_2 q_2 g(\dfrac{x-x_{1}}{q_2}) + y_{1} (x-x_1)  +p_1 q_1 g(1) \end{cases}$, \\where $g(1)=\dfrac{y_{1} q_2}{1- p_1 q_1 -p_2 q_2 }$.

It seems that we can derive the points $y_1$, $y_2$ from $g_i=g(x_i)$:

$y_1=\left(\dfrac{1}{q_1}-\dfrac{1}{q_2}\right)g_1 + \left(\dfrac{1}{q_2} - 1\right)g_2$, $y_2=\dfrac{q_1 g_2}{g_1} y_1$.

So, to construct a zipper with the vertices $(0,0),\ (x_1,g_1),\ (1,g_2)$ and the signature $\e=(0,0)$ with smooth attractor, we can take

$W_1 (x,y)=\left( \dfrac{x}{q_1} , \dfrac{g_1}{g_2}y \right)$, 

$W_2 (x,y)=\left( \dfrac{x}{q_2} + x_1 , \left(1-\dfrac{g_1}{q_1 g_2}\right)q_2 y + \left(\left(\dfrac{1}{q_1}-\dfrac{1}{q_2}\right)g_1 + \left(\dfrac{1}{q_2} - 1)g_2\right) \right) \left( \dfrac{x-x_1}{q_2} - x_1\right) +g_1\right)$.
\\

{\bf Example 2.} Consider another case of zipper on the plane, consisting of:

\begin {center}
$\begin {cases}
S_{1}  \left( \begin{array}{c} x \\ y \\ \end{array} \right) = p A \left( \begin{array}{c} x \\ y \\ \end{array} \right) \\
S_{2}  \left( \begin{array}{c} x \\ y \\ \end{array} \right) = p B \left(\ \begin{array}{c} x \\ y \\ \end{array}
\right) + \left(\ \begin{array}{c} \dfrac{1}{2} \\ h \\ \end{array} \right) \end{cases}$
\end {center}
where $ 0<h<\sqrt{3}/2$, $p = \sqrt{h^{2} + \dfrac{1}{4}}$, $\alpha = \arctan (2h)$, $A$ and $B$ are the rotation matrixes:

$A = \left( \begin{array}{cc} \cos \alpha & - \sin \alpha \\ \sin \alpha & \cos \alpha\\ \end{array} \right) $,
$B = \left( \begin{array}{cc} \cos \alpha & \sin \alpha \\ - \sin \alpha & \cos \alpha\\ \end{array} \right) $.

However, $S_{1}: [0,1] \times [0,0] \rightarrow \left[0, \dfrac{1}{2}\right] \times [0,h]$, $S: [0,1] \times [0,0] \rightarrow \left[\dfrac{1}{2},1\right] \times [0,h]$.

Parametrize zipper so that $f\left( \dfrac{1}{2} \right)=\left(\begin{array}{c} 1/2 \\ h \end{array}\right)$:

\begin {center}
$f(t) = \begin {cases} pAf(2t) , t\in \left[0;\dfrac{1}{2}\right] \\
pBf(2t-1) + \left(\ \begin{array}{c} \dfrac{1}{2} \\ h \\ \end{array} \right) , t\in \left[\dfrac{1}{2};1\right] \end {cases}$
\end {center}

Let $g(t) = \int\limits_{0}^{t} f(\xi)d\xi$.

Then:

\begin {center}
$g(t) = \begin {cases} \dfrac{p}{2}Ag(2t), t\in \left[0;\dfrac{1}{2}\right]\\
\dfrac{p}{2}B g(2t-1) + \dfrac{p}{2} A \left(\ \begin{array}{c} \dfrac{1}{2(1-p\cos\alpha)} \\ \dfrac{h}{1-p\cos\alpha} \\
\end{array} \right) + \left(\ \begin{array}{c} \dfrac{t}{2} \\ ht \\ \end{array}
\right) , t\in \left[\dfrac{1}{2};1\right]
\end {cases}$
\end {center}

The corresponding zipper consists of two mappings that transform unit cube into a rectangular parallelepiped
$\left[0,\dfrac{1}{2}\right]\times\left[0,\dfrac{1}{2}\right]\times[0,h]$ and $\left[\dfrac{1}{2},1\right]\times\left[\dfrac{1}{2},1\right]\times[0,h]$:

\begin {center}
$\begin {cases}
W_{1}\left( \begin{array}{c} t \\ x \\ y \\ \end{array} \right) = \left(  \begin{array}{c} \dfrac{t}{2} \\ \dfrac{p}{2}A \left(\ \begin{array}{c} x \\ y \\ \end{array} \right)  \end{array} \right)\\
W_{2}\left( \begin{array}{c} t \\ x \\ y \\ \end{array} \right) =  \left( \begin{array}{c} \dfrac{1}{2} \\ \dfrac{t}{2} \\ ht  \end{array} \right) + \left( \begin{array}{c} \dfrac{t}{2}  \\ \dfrac{p}{2}B \left(\ \begin{array}{c} x \\ y \\ \end{array} \right)  \end{array} \right)
\end {cases}$
\end {center}


\newpage


\begin{thebibliography}{}


\bibitem{1}  V.V. Aseev, A.V. Tetenov, A.S.Kravchenko  On  Self-Similar Jordan Arcs in Plane / 
Siberian Mathematical Journal
 2003, Volume 44, Issue 3, pp 379-386 

\bibitem{2} V.V. Aseev, A.V. Tetenov, On the Self-Similar Jordan Arcs Admitting Structure Parametrization / 
Siberian Mathematical Journal
 2005, Volume 46, Issue 4, pp 581-592 


\bibitem{3}    M.F. Barnsley, Fractal functions and interpolation / Constructive Approximation, 1986, 2, 303-329.


\bibitem{4}  J. Hutchinson, Fractals and self-similarity / Indiana University Mathematics Journal, 1981, 30, 713-747.

 \bibitem{5}   A.S.Kravchenko , Smooth self-affine  zippers. Sobolev Math Institute  preprint, Novosibirsk, 2005.(Russian)


\end{thebibliography}
\end{document}